%% file: MCM2014_arxiv_version.tex
\documentclass{article} % For LaTeX2e

\input{definitions_arxiv.tex}

\title{What Makes a Good Plan?\\An Efficient Planning Approach to Control Diffusion Processes in Networks}

\newcommand{\myaffiliation}{CMLA {--} ENS Cachan, CNRS, France}

\author{
Kevin Scaman
\hspace{2em}
Argyris Kalogeratos
\hspace{2em}
Nicolas Vayatis
\\\myaffiliation\\
\texttt{\{scaman, kalogeratos, vayatis\}@cmla.ens-cachan.fr}
}

\begin{document}
%=================================================================
%
\maketitle
%
%
%
%
%=================================================================
\begin{abstract}
%=================================================================
%

In this paper, we analyze the quality of a large class of simple dynamic resource allocation (DRA) strategies which we name \emph{priority planning}. Their aim is to control an undesired diffusion process by distributing resources to the contagious nodes of the network according to a predefined \emph{priority-order}. In our analysis, we reduce the DRA problem to the \emph{linear arrangement} of the nodes of the network. Under this perspective, we shed light on the role of a fundamental characteristic of this arrangement, the \emph{maximum cutwidth}, for assessing the quality of any priority planning strategy. Our theoretical analysis validates the role of the maximum cutwidth by deriving bounds for the extinction time of the diffusion process. Finally, using the results of our analysis, we propose a novel and efficient DRA strategy, called \emph{Maximum Cutwidth Minimization}, that outperforms other competing strategies in our simulations.

\end{abstract}

%
%
%
%
%=================================================================
\section{Introduction}\label{sec:intro}
%=================================================================
%
Diffusion processes usually arise in systems involving agents whose behaviors depend on their close environments. Diseases, computer viruses, ideas and interests are spread through a social network by means of interactions among its users. In all these phenomena, a change in one agent may affect the actions of other neighboring agents, resulting, under certain conditions, in a massive change of behavior at the network scale.

A fair body of articles consider the problem of influence maximization, which attempts to maximize the spread of a diffusion process. On the other hand, being able to suppress or remove an undesired information or social diffusion process has received less attention, even though critical in many real-life situations. In public health, scenarios in which the spread of a virus needs to be controlled have been extensively studied by epidemiologists. Various analogues emerged in modern information networks in which a diffusion can be engineered to be viral and may cause huge negative social and economic effects. For instance, during the London riots during the summer of 2011, public opinion was obfuscated and baffled by many false rumors \cite{UKRIOTS2013}. During large scale natural disasters, there is also a risk to follow misinformation diffused through social media by various individuals, while trying to coordinate some rescue and volunteer teams \cite{MDTNHE2010, OPCPFSDHS2014}. 
Another recent example arose in March 2013, when a false tweet from the account of Associated Press was headlined with ``\emph{Breaking: Two explosions in the White House and Barack Obama is injured}". The false rumor was immediately retweeted by about 2 million users causing panic and speculation for a few minutes in the US stock market. As a matter of fact, negative publicity is highly damaging for organizations and brands: even though ill-founded, a slander can dramatically affect the interests of a company due to the massive scale of the buzz-like effect. 

%-- related work
The control of diffusion processes has been studied in various fields, including epidemiology and computer networks resilience; the respective literature can generally be divided in three complementary lines of research:
\begin{itemize}
\item \emph{Static vaccination strategies}. Most of the epidemic literature focuses on static control actions such as permanently removing a set of edges or nodes of the network \cite{cohen2003efficient, tong2012gelling, wang2003epidemic, schneider2011suppressing, preciado2013optimal}. In this case, the available budget is considered fixed, and the effect of a control action \emph{permanent}. Examples of static resource allocation strategies can be found in \cite{preciadoZEJP13, pageRankAntidotes2009}.
\item \emph{Budget optimization}. Complementary to resource allocation, the determination of the right budget of resources to be spent each time step, in order to fulfill cost and efficiency constraints, has critical impact on the resulting strategy. 
Several studies lay on the line of optimal budget estimation, assuming that the network administrator is capable of storing resources for later use \cite{klepac2012, forster2007optimizing, khouzani2011optimal}. These approaches usually make the simplifying assumption of \emph{uniform mixing}, \ie that the contagious nodes are uniformly scattered in the network. Therefore, they do not address the problem of how exactly to allocate the resources on the nodes of the network, but rather how many resources should be provided at each time to cause a desired macroscopic result.
\item \emph{Dynamic resource allocation}. A few studies consider dynamic strategies for allocating resources for dealing with epidemics. One of the most well-known such strategy is \emph{contact-tracing} \cite{RSA:RSA20315} that consists in healing the neighbors of contagious nodes. In practice, this approach has been shown inefficient to contain epidemics \cite{RSA:RSA20315}, especially when they are beyond a very initial state.
\end{itemize}

Our major contribution in this article is to introduce a particular class of strategies for suppressing an undesired diffusion. Instead of choosing a set of nodes whose behavior will be permanently modified, we allow the network administrator to change the distribution of the \emph{resources} during the diffusion. In other words, we consider \emph{targeted} and \emph{temporal} action on individual nodes of the network, that can affect their behavior. Since reacting to fast spreading phenomena is difficult to achieve, we consider a simple set of dynamic strategies relying on a predefined \emph{priority-order}. The strategy will gradually suppress the diffusion and finally remove the undesired contagion by focusing on the first contagious users according to the priority order.

In what follows, \Sec{sec:model} describes the model used for the diffusion process as well as the control actions available to the network administrator. \Sec{sec:priorityplanning} presents the idea of \emph{priority planning} as well as a natural representation of the problem connecting our analysis to \emph{linear arrangement} problems. This anaylysis sheds light on the role of the network's \emph{maximum cutwidth} for the efficiency of a given strategy. By minimizing this value, we develop an efficient strategy in \Sec{sec:algorithm}, and validate in \Sec{sec:theory} that our intuition is valid by deriving theoretical results on the extinction time of the diffusion process. 
Finally, we present experimental results in \Sec{sec:exps} and show that: i)~The derived bounds are very close to the epidemic threshold, thus validating the fundamental role of the \emph{maximum cutwidth} in the evaluation of such strategies. ii)~The proposed strategy outperforms its competitors in a wide range of scenarios.

%
%
%=================================================================
\section{Diffusion and control model}\label{sec:model}
%=================================================================
%
%--------------------------------------------
\subsection{The \emph{Susceptible-Infected-Susceptible} epidemic model}
%--------------------------------------------
%
We will consider a simple diffusion process known in the epidemiology literature as \emph{N-Intertwined Susceptible-Infected-Susceptible} (SIS) model \cite{van2009virus}. According to this model, a diffusion can be spread through the edges of the network and turn the state of nodes from \emph{susceptible} (or \emph{healthy}) to \emph{infected}. When a node becomes infected, it can in turn spread the contagion to its direct neighbors and, after some amount of time, return to the susceptible state without bearing any immunity. At each time, the network administrator can take control actions in order to reduce the epidemic. These actions are represented as a \emph{budget} of \emph{resources} (or treatments in the epidemic analog) to distribute in the network. Each resource increases the recovery rate of the receiver node. In a continuous-time framework, the distribution of control resources can be revised anytime. However, in practice the situations in which this is needed is only when there is a change in the state of the network, \ie a new node infection or recovery, or a modification in the available resource budget. 

This model is well suited to situations in which an undesired contagion affects a network, and the control action is \emph{local} and \emph{expensive}. Among other application examples, controlling epidemics using antidotes, limiting rumors via targeted action or allocating resources geographically to fight against a societal problem, seem valid scenarios for such a diffusion and control model.

\spinlinetitle{Formal definition}{.} 
Let $\mathcal{G} \op{=} (\mathcal{V}, \mathcal{E})$ be a network of $N$ nodes with adjacency matrix $A$, where $A_{ij} \op{=} 1$ if $i \op{\neq} j$ and edge $(i,j) \op{\in} \mathcal{E}$, else $A_{ij} \op{=} 0$. We describe the state of the diffusion process with a \emph{state vector} $X(t) \op{\in} \real^N$ that keeps the state of each node of the network: $X_i(t) \op{=} 1$ if node $i$ is contagious at time $t$, else $X_i(t) \op{=} 0$. We also describe the control action on the network with a \emph{resource vector} $R(t)$, where $R_i(t)$ is $1$ iff node $i$ is given a resource at time $t$. In such case, and following the epidemic analogy, we say that node $i$ is being \emph{healed} by the resource. Using also the formalism of \cite{ganesh2005effect}, we model the diffusion with a \emph{continuous-time Markov process} which has the transition rates:
\begin{equation} \label{eq:Xt} \formulastyle
\begin{array}{l}
X_i(t) : 0\rightarrow1 \mbox{ at rate } \beta \sum_j {A_{ji} X_j(t)};\\
X_i(t) : 1\rightarrow0 \mbox{ at rate } \delta + \rho R_i(t),\\
\end{array}
\end{equation}
where $\beta$, $\delta$, $\rho$ are, respectively, the transmission rate over one network edge, the recovery rate without receiving a resource, and the increase in recovery rate that a resource induces. Note that the diffusion process is continuous in time and thus $X(t)$ is a stochastic process. The action of the resource on a node increases its chances to return to the \emph{healthy} state. Finally, we define two dimensionless parameters: $r \op{=} \frac{\beta}{\delta}$ the \emph{effective spreading rate} of the diffusion, and $e \op{=} \frac{\rho}{\delta}$ the \emph{resource efficiency}.

A \emph{control strategy} takes as input the network $\mathcal{G}$, the characteristics of the diffusion ($\beta$, $\delta$, and $\rho$) and the network state $X(t)$, and returns the distribution of the resources $R(t)$ (\Eq{eq:Xt}). In the following, we refer to the problem of finding an optimal control strategy with respect to the minimization of the spread of the diffusion process as the \emph{dynamic resource allocation} (DRA) problem. We consider as \emph{budget} $b(t)$ the maximum number of resources that can be distributed in the network at time $t$, and that the available budget at time $t$ cannot be stored for later use. Following the epidemiology literature, we will denote as \emph{epidemic threshold under a given strategy} the resource efficiency $e$ above which the control actions removes the contagion in \emph{reasonable time}, that is less than exponential in the number of nodes of the network.
%
%--------------------------------------------
\subsection{Priority planning: a control plan to gradually remove a contagion}\label{sec:priorityplanning}
%--------------------------------------------
%
%------------------------------------------------------
\subsubsection{Definition}
%------------------------------------------------------
%
A \emph{priority planning} is a DRA strategy that distributes resources to the top-$b(t)$ infected nodes according to a fixed \emph{priority ordering} of the nodes, where $b(t)$ is the available budget of resources at each time. In mathematical terms, an \emph{$\lorder$-priority planning} is defined by a bijective mapping $\lorder\!: \mathcal{V} \op{\rightarrow} \{1, ..., N\}$ of the $N$ nodes of the network \st $\lorder(v)$ is the position of node $v$ in the priority order. The strategy then selects the first $b(t)$ infected nodes according to $\lorder$:
\begin{equation}
R_i(t) = \left\{
\begin{array}{ll}
1 &\mbox{if } X_i(t) = 1 \mbox{ and } \lorder(v_i) \leq \theta(t);\\
0 &\mbox{otherwise},\\
\end{array}\right.
\end{equation}
where $\theta(t)$ is an adjusted threshold set so that $\sum_i R_i(t) \op{=} b(t)$.

These strategies may be regarded as simple planning strategies for the removal of an undesired contagion: the order in which the nodes will be healed is determined prior to the beginning of the diffusion. Then, this order is respected no matter how the diffusion process evolves, and the strategy removes the contagious nodes gradually starting from the first in the priority-order to the last. Although our strategy (\Sec{sec:algorithm}) can handle a generic time-dependent budget, we only consider a constant budget rate $b(t) \op{=} b_{tot}$ for theoretical developments as well as focused experimental evaluation.

%
%------------------------------------------------------
\subsubsection{Interpretation of priority planning}\label{sec:interpretation}
%------------------------------------------------------
%
In this section, we present a novel perspective on analyzing the DRA problem which leads to efficient priority planning strategies. More specifically, we reduce the problem of determining a good priority-order to a \emph{linear arrangement} (LA) problem. 

\begin{figure}[!t] \footnotesize
\beforecaptvskip
  \begin{subfigure}[b]{0.5\textwidth}
		\centering
		\includegraphics[viewport=20 545 330 636, clip=true, width=0.92\linewidth]{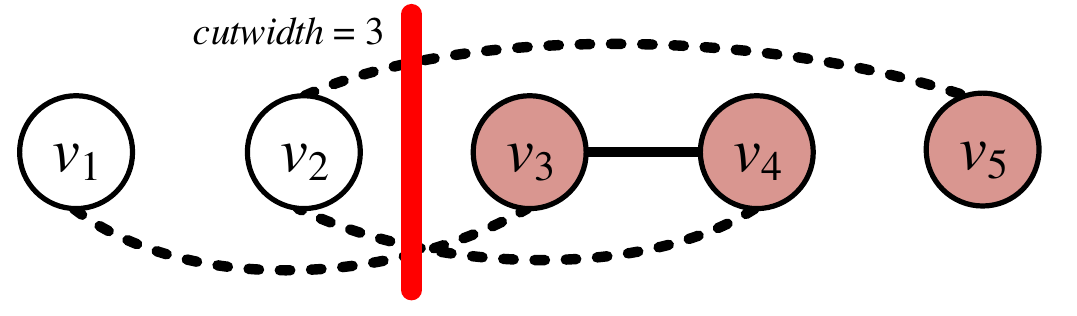}
		\caption{LA with mapping $\lorder\!: \mathcal{V} \op{\rightarrow} \{1,2,3,4,5\}$}
		\label{fig:badLA}
	\end{subfigure}
	\begin{subfigure}[b]{0.5\textwidth} 
		\centering
		\includegraphics[viewport=20 640 330 734, clip=true, width=0.92\linewidth]{network_example2.pdf}
		\caption{LA with mapping $\lorder'\!: \mathcal{V} \op{\rightarrow} \{1,3,4,2,5\}$}
		\label{fig:goodLA}
	\end{subfigure}	
	\beforecaptvskip
	\caption{The priority orders of two static strategies are visualized as the linear arrangements of the nodes based on two mappings $\lorder$ and $\lorder'$. Healthy nodes (white) and contagious nodes (red) form two distinct subgraphs that are connected by the \emph{contagious edges} (dashed lines) of the \emph{front} illustrated as a red vertical line. (a)~M$p$LA cost, for $p$ \op{=} 1, $\objf_1(\mathcal{G}, \lorder) \op{=} 8$ and MCLA cost $\objf_2(\mathcal{G}, \lorder) \op{=} 3$; (b)~$\objf_1(\mathcal{G}, \lorder') \op{=} 4$ and $\objf_2(\mathcal{G}, \lorder') \op{=} 1$, respectively.}
	\label{fig:static_strategy}
\end{figure}

\spinlinetitle{Linear arrangement}{.}
Formally, a linear arrangement $\lorder\!: \mathcal{V} \op{\rightarrow} \{1, ..., N\}$ maps the nodes of $\mathcal{G}$ on $N$ discrete positions located on a line by assigning one position-label to each node (\Fig{fig:static_strategy}). LA is a class of combinatorial optimization problems whose purpose is to minimize some functional $\objf$ over the space $\mathcal{L}$ of all possible node permutations: $\lorder^* \op{=} \argmax_{\lorder \in \mathcal{L}} \objf(\mathcal{G},\lorder)$. The problems in the LA class are also referred to as \emph{graph layout} problems (for a survey see \cite{SGLP2002}), and indicative applications are the graph drawing, VLSI design, and network scheduling.

Probably the most popular LA instance is the \emph{minimum $p$-sum linear arrangement} (M$p$LA, or MLA when $p \op{=} 1$) \cite{OANV1964} that minimizes the following functional:
\begin{equation} \label{eq:p-sum}
\begin{array}{l}
\mbox{M$p$LA}\!: \ \ \objf_1(\mathcal{G}, \lorder) = \Big(\sum_{i,j} A_{ij}\left|\lorder(v_i) - \lorder(v_j)\right|^p\Big)^{1/p}.
\end{array}
\end{equation}
In other words, for $p \op{=} 1$, this minimizes the weighted sum over all distances between node assignements $\ell(u)$ to the $N$ discrete equally spaced positions. Another category of LA problems is related to the minimization of some maximum value of the LA, for example the \emph{bandwidth} (\ie the maximum edge length), or the \emph{workbound} (\ie the sum of maximum edge length of each node). There is also the \emph{minimum cutwidth linear arrangement} (MCLA) problem, where the \emph{cutwidth} at a location $c$ lays among the node positions $c$ and $c \op{+} 1$ of the LA and is the weighted edge-cut at that position. Then, the minimized functional is the maximum of the $N \op{-} 1$ local cuts, and more precisely:
\begin{equation} \label{eq:cw}
\begin{array}{l}
\mbox{MCLA}\!: \ \ \objf_2(\mathcal{G}, \lorder) = \max_{c={1,...,N-1}} \sum_{i,j} A_{ij} \one_{\{\lorder(v_i) \leq c < \lorder(v_j)\}},
\end{array}
\end{equation}
where $\one_{\{\cdot\}}$ is the indicator function that returns one if the input condition is true, else zero. For the simplicity of notations we also denote the value of the above functional as $\MaxCut \op{=} \objf_2(\mathcal{G}, \lorder)$. Note that symmetric LAs are evaluated as of equal quality by $\objf_1$, $\objf_2$.
\Fig{fig:static_strategy} illustrates two LAs of a small network and an example of the cost values derived in each case by \Eq{eq:p-sum} and \Eq{eq:cw}. Among them, the arrangement of \Fig{fig:goodLA} is better, indeed optimal, in terms of both objective functions $\objf_1$, $\objf_2$.

\spinlinetitle{Designing an efficient control plan}{.}
The analysis of the diffusion process under the perspective of a linear arrangement, sheds light on how to design a DRA control plan to suppress a diffusion efficiently. More specifically, knowing that the contagion can only be spread through the \emph{contagious edges} which connect infected with healthy nodes, and that any DRA plan is based on a static priority-order for healing the nodes, we realize that it is highly critical to have as less contagious edges as possible during the whole diffusion process.

\Fig{fig:static_strategy} represents two diffusion instances that are being suppressed according to different plans; we consider that each priority-order is deployed from the left to the right. The current state of the plan is denoted with a red vertical line, the \emph{front}, that roughly separates the healthy from the contagious nodes and indicates where the strategy would allocate resources. In practice, healthy nodes can also appear on the right side of the front due to self-recovery. The number of edges crossing through the front, \ie the cutwidth at that position, indicates how vulnerable are actually the healthy nodes of the cleared part of the network laying left to the front. Therefore, we can argue that a good solution to the M$p$LA minimization would provide a smoother plan that is generally \emph{easier to proceed from state-to-state}, while the MCLA would minimize the $\MaxCut$ which is \emph{the most difficult state of the plan} and, thus, a determinant for its eventual success in completely removing the contagion. This implies that the minimization of the $\MaxCut$ can be a proxy for minimizing the epidemic threshold for any arbitrary network. Next, we theoretically show that these intuitive remarks are valid. 
%
%
%
%=================================================================
\section{The \emph{Maximum Cutwidth Minimization} control strategy}\label{sec:algorithm}
%=================================================================
%
Based on the analysis of the previous section, that uncovered a strong dependency between the epidemic threshold of the diffusion process and the \emph{maximum cutwidth}, we propose a new DRA algorithm for arbitrary networks where the main idea is to distribute resources to contagious nodes in the order that minimizes $\MaxCut$. Specifically, given a network $\mathcal{G}$, we compute, prior to the diffusion process, the linear arrangement $\lorder_{MC}(\mathcal{G})$ of its nodes with minimum $\MaxCut$ value using any available optimization algorithm for this problem. Then, during the diffusion, the strategy distributes the budget of resources to the contagious nodes according to the order of $\lorder_{MC}(\mathcal{G})$. The pseudocode of the strategy is summarized in \Alg{alg:maxcutmin}, while theoretical justification and results are provided in \Sec{sec:theory}.

\begin{algorithm}[tb!] \footnotesize
   \caption{\emph{Maximum cutwidth minimization} (MCM) control strategy}
   \label{alg:maxcutmin}
\begin{algorithmic} \small
	 \STATE \textbf{Input\:\:\ :}\hspace{0.4em}network $\mathcal{G}$, infection state vector $X(t)$, budget size $b(t)$
	 \STATE \textbf{Output:}\hspace{0.4em}the resource allocation vector $R(t)$
	 \vspace{0.6em}
	 \STATE \raisebox{0.1em}{$\rhd$}~\textbf{Prior to the diffusion:}
	 \STATE Compute the priority-order $\lorder \op{=} \lorder_{MC}(\mathcal{G})$ by minimizing the maximum cutwidth $\MaxCut$
	 \STATE Order the nodes of $\mathcal{G}$ according to $\lorder$, \ie compute the node list $(v_1,...,v_N)$ \st $\forall i \op{\in} \{1,..., N\}$, $\lorder(v_i) \op{=} i$ \\
	 \vspace{0.7em}
	 \STATE \raisebox{0.1em}{$\rhd$}~\textbf{During the diffusion:}
	 \IF{$\sum_i X_i(t) < b(t)$}
	 \RETURN $X(t)$
	 \ENDIF
	
	 \STATE Let $R(t)\op{=}\zero$ a zero $N$-dimensional vector, also $i \op{=} 1$ and $budget \op{=} b(t)$
	 \WHILE{$budget > 0$}
	 \IF{$X_{v_i}(t) = 1$}
	 \STATE $R_{v_i}(t) \leftarrow 1$
	 \STATE $budget \leftarrow budget - 1$
	 \ENDIF
	 \STATE $i \leftarrow i + 1$
	 \ENDWHILE
	
	 \RETURN $R(t)$
\end{algorithmic}
\end{algorithm}

\spinlinetitle{Maximum cutwidth linear arrangement}{.}
MCLA is known to be an NP-hard problem. However, approximation heuristics do exist in the related literature \cite{SSCMP2012, VFSCMP2013}. One of the difficulties of this problem is that the cost function to optimize (\Eq{eq:cw}) is extremely flat in the search space, due to the fact that slight changes in the arrangement will most probably not change $\MaxCut$. For this reason, we chose to relax the MCLA problem by optimizing the \emph{p-sum linear arrangement} problem with $p \op{=} 1$ (MLA, see \Sec{sec:interpretation}), which is easier than MCLA and more suited to gradient descent or simulated annealing methods. Furthermore, as discussed in \Sec{sec:interpretation}, MLA may produce a smooth priority-order that exhibits some desirable properties. 

\spinlinetitle{Scalability of the method}{.}
The scalability of the strategy is highly dependent on the chosen algorithm for finding the optimal order. For our experimental results, we applied a hierarchical approach to take advantage of the group-structure of the social network: i)~first, we identified dense clusters by applying \emph{spectral clustering} and we ordered the clusters (considered as high-level nodes) using \emph{spectral sequencing} \cite{OLLEG1992}, ii)~then, we computed a good ordering of the nodes in each of the clusters independently using spectral sequencing followed by an iterative approach based on random node swaps, inspired by \cite{ETSSAAMLA2008}, iii)~finally, the swap-based approach was reapplied to optimize the ordering all together. The whole process achieved fairly good results (see \Tab{tab:twittermaxcuts}) in reasonable time. Since clustering and spectral sequencing depend on the computation of eigenvectors for the highest eigenvalues of an $N \op{\times} N$ sparse matrix with $|\mathcal{E}|$ non-zero entries, the overall algorithm has a complexity $O(|\mathcal{V}| \op{+} |\mathcal{E}|)$ \cite{Arora05fastalgorithms}. Hence, MCM is generally scalable to the size of real social networks.

%
%
%=================================================================
\section{Theoretical bounds for the extinction time and the epidemic threshold}\label{sec:theory}
%=================================================================
%
We now provide a justification of the design of our MCM algorithm. The following theorem gives an \emph{upper bound} for the expected extinction time under any priority planning in the simple case of $b(t) \op{=} 1$. The general case, while much more complex to analyze, is very similar in our experiments provided that the treatment efficiency $e$ is multiplied by the available budget $b(t)$ (see \Sec{sec:exps}). %# healed or treated? 
Above a threshold value that depends on the \emph{maximum cutwidth} $\MaxCut$ of the network under a considered priority-order, \Theorem{th:boundExtinctionTime} proves that the diffusion process converges in reasonable time to its absorbent state.

\begin{theorem}\label{th:boundExtinctionTime}
Let $\mathcal{G} \op{=} (\mathcal{V}, \mathcal{E})$ be a network of $N$ nodes and assume a fixed budget $b(t) \op{=} 1$. Let $\lorder\!: \mathcal{V} \op{\rightarrow} \{1,..., N\}$ be an ordering of the nodes of $\mathcal{G}$, and $\MaxCut$ be the maximum cutwidth of $\lorder$, or equivalently the highest number of contagious edges during the planned removal of the contagion. Then the following upper bound holds for the expected extinction time $\Exp{\tau}$ under the $\lorder$-priority plan and starting from a total infection:

If $\rho > \beta \MaxCut \left(1 +  2\sqrt{\epsilon} + \epsilon \right) - \delta$, then
\\
\vspace{-0.3em}
\begin{equation}
\Exp{\tau} \ \leq \ \frac{N}{\rho + \delta - \beta \MaxCut \left(1 + 2\sqrt{\epsilon} + \epsilon \right)}
\end{equation}
where $\epsilon = \frac{d_{max} (1 + \ln{N})}{\MaxCut}$ and $d_{max}$ is the maximum node degree of the network.
\end{theorem}
This bound relates the extinction time to the number of contagious edges in the worst step of the strategy (\ie its \emph{maximum cutwidth}). When $\MaxCut$ is such that $d_{max} (1 \op{+} \ln{N}) \op{\ll} \MaxCut$, this formula bounds the epidemic threshold under a given priority planning by $r \MaxCut \op{-} 1$, where $r \op{=} \beta / \delta$ is the \emph{effective spreading rate}, thus giving a special emphasis on the role of $\MaxCut$ for assessing the quality of a given strategy.
\begin{corollary}\label{th:boundThreshold}
Under the same hypotheses as \Theorem{th:boundExtinctionTime}, the epidemic threshold $e^*$, \ie the lowest resource efficiency $e \op{=} \frac{\rho}{\delta}$ \st the diffusion converges to zero in reasonable time, is upper bounded:
\begin{equation}
e^* \ \leq \ r \MaxCut \left(1 + 2\sqrt{\epsilon} + \epsilon \right) - 1,
\end{equation}
where $\epsilon \op{=} \frac{d_{max} (1 + \ln{N})}{\MaxCut}$ and $r \op{=} \frac{\beta}{\delta}$ is the effective spreading rate.
\end{corollary}

While $\epsilon$ is necessary for the theorem to hold, we rather consider it as an artifact of the proof and not a fundamental aspect of the result. As such, we expect $r \MaxCut \op{-} 1$ to be closer to the epidemic threshold (see specifically \Sec{sec:qualityExps} for an experimental validation of this intuition). This result verifies that, according to the intuition provided in \Sec{sec:interpretation}, removing an undesired contagion requires the resource efficiency to be \emph{as high as needed for the worst step} of the specified plan. Equivalent to the celebrated relationship between the epidemic threshold under no control strategy and the spectral radius of the adjacency matrix, this result is fundamental for understanding the behavior of the diffusion process and designing efficient DRA strategies. The simulations in \Sec{sec:exps}, attest that minimizing this upper bound is very efficient for dynamically controlling a diffusion process.
%
%
%
%
%=================================================================
\section{Experimental results}\label{sec:exps}
%=================================================================
%
%--------------------------------------------
\subsection{Quality of the theoretical bound}\label{sec:qualityExps}
%--------------------------------------------
%

\begin{figure}[t] \footnotesize
\beforecaptvskip
	\centering
	\begin{subfigure}[b]{0.45\textwidth}
		\centering 
		\includegraphics[viewport=17 4 381 306, clip=true, width=1\linewidth]{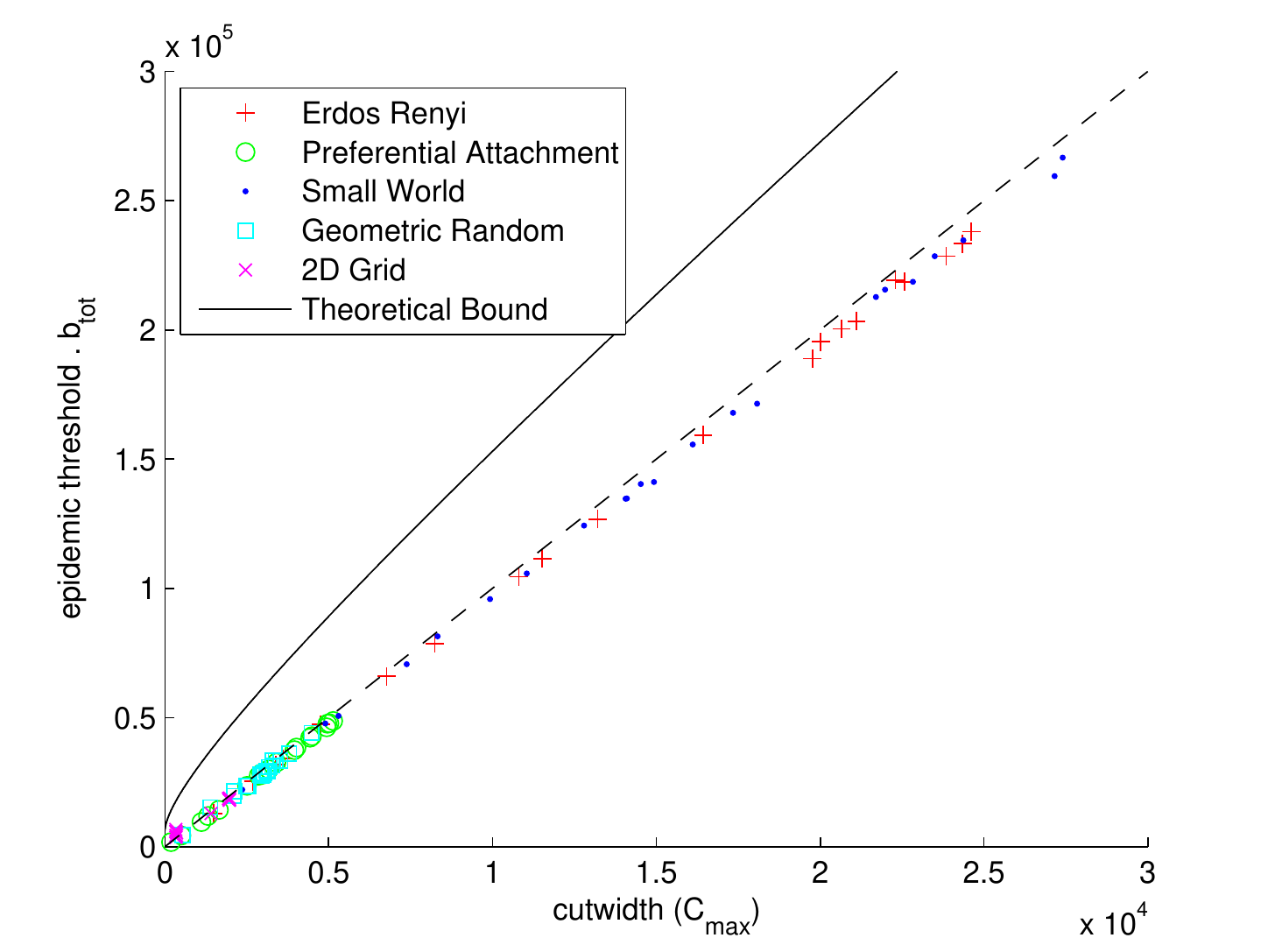}
		\caption{high infectivity scenario: $r \op{=} 10$}
	\end{subfigure}%
	\hspace{2em}
	\begin{subfigure}[b]{0.45\textwidth}
		\centering
		\includegraphics[viewport=9 4 381 306, clip=true, width=1\linewidth]{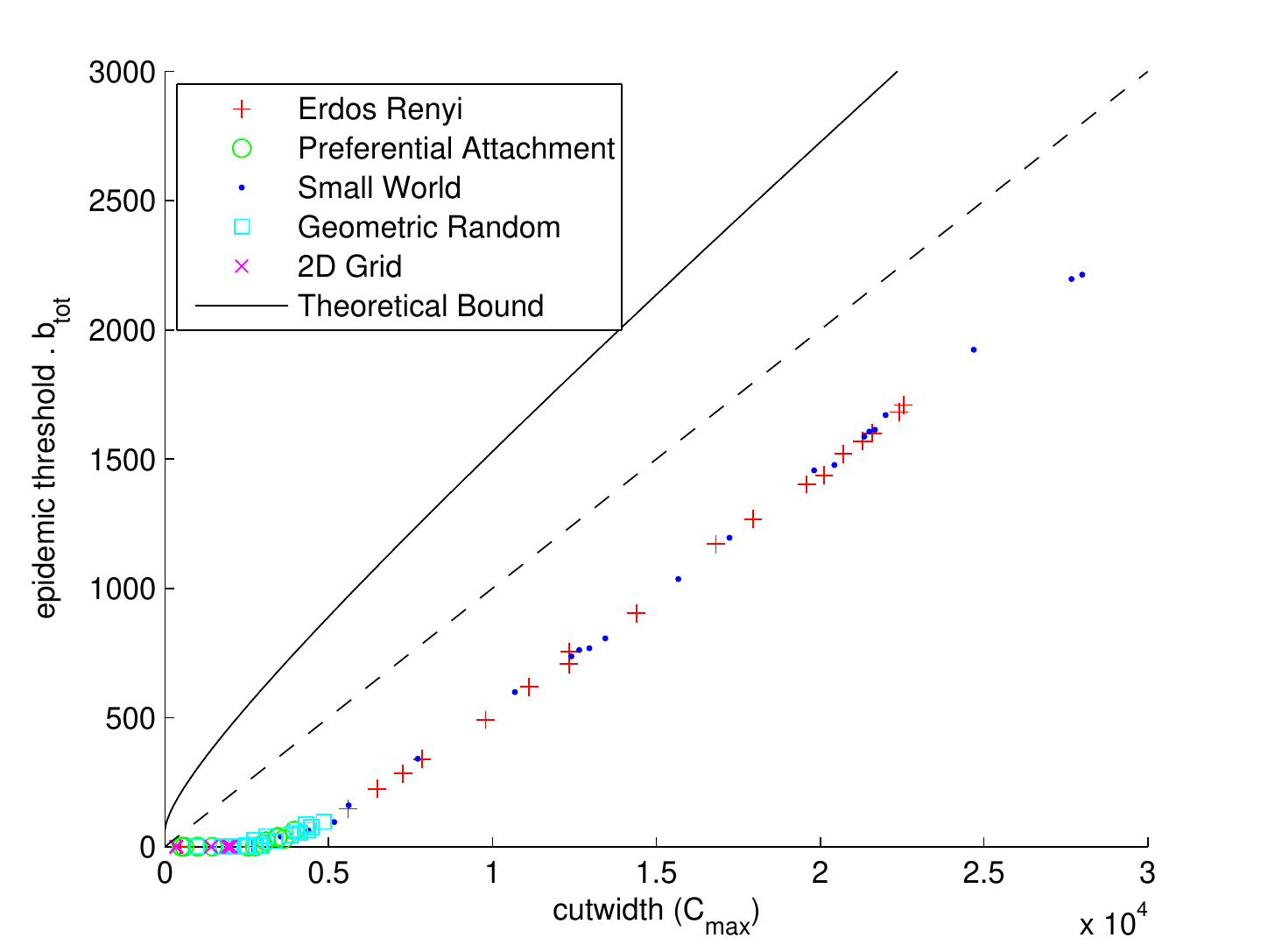}
		\caption{low infectivity scenario: $r \op{=} 0.1$}
	\end{subfigure}
	\beforecaptvskip
	\caption{Epidemic threshold \wrt the maximum cutwidth $\MaxCut$ for various random network types of $N \op{=} 1000$ nodes. For the theoretical bound, $d_{max} \op{=} 100$ was used. The dashed line indicator has slope equal to the effective spreading rate $r$.}
	\label{fig:boundQuality}
\end{figure}

In \Fig{fig:boundQuality}, we show that the relationship between the \emph{maximum cutwidth} $\MaxCut$ and the epidemic threshold under a specific priority plan is very stable, nearly linear, and hence $\frac{r\MaxCut}{b_{tot}}$ is in fact a very good approximation. Each plotted point is the epidemic threshold under a given strategy and for a given network instance. We sampled $100$ networks of $1000$ nodes from 5 random network generators (\eg see \cite{Newman:2010:NI} for details on these network types): i)~\Erdos, ii)~Preferential attachment, iii)~Small-world, iv)~Geometric random \cite{penrose2003random} and v)~2D regular grids. We also sampled DRA control strategies at random among the set of strategies described next in \Sec{sec:twitterexps}, in order to cover a wide range of different scenarios.

\Fig{fig:boundQuality} is a summary of these results. The epidemic threshold is always below $\frac{r \MaxCut}{b_{tot}}$, but it stays very close to this value and, more importantly, shows a close to deterministic behavior with respect to this value even in the case of low infectivity (small $r$ value) where intuitively the random self-recovery of nodes become more significant. This result justifies the minimization of $\MaxCut$ as a proxy for minimizing the epidemic threshold, and thus remove the contagion more efficiently.
%
%
%
%--------------------------------------------
\subsection{Application on a real social network}\label{sec:twitterexps}
%--------------------------------------------
%
%
%
\begin{table}[!t] \footnotesize
\center
\begin{tabular}{l|l|c|c} 
\hline
\textbf{Strategy} & \textbf{Max cutwidth} & \textbf{Max cutwidth} & \textbf{Expected epidemic threshold}\\
 \textbf{} & \ \ \ \ \ \ \ \ \ \  & \!\textbf{\% MCM} & \textbf{($r \op{=} 0.1$, $b_{tot} \op{=} 100$)}\\
\hline \hline
	RAND & 670,000 $\pm$ 1000 & 931\% & 670\\
	MN     & 628,571            & 874\% & 629\\
	LN     & 628,571            & 874\% & 629\\
	LRSR   & 349,440            & 486\% & 349\\
	\textbf{MCM}    & \ \ \textbf{71,956}            & \textbf{100\%} & \ \ \textbf{72}\\
\hline
\end{tabular}
\beforecaptvskip
\caption{Maximum cutwidth values ($\MaxCut$) for different strategies in the TwitterNet used in our experiments. The expected epidemic threshold is $\frac{r\mathcal{C}_{max}}{b_{tot}}$ based on \Corollary{th:boundThreshold} and the experiments of \Sec{sec:qualityExps}.}
\label{tab:twittermaxcuts}
\end{table}%
\begin{figure}[t]
\vspace{-0.2em}
	\centering
	\begin{subfigure}[b]{0.45\textwidth}  
		\centering
	\includegraphics[viewport=12 0 405 307, clip=true, width=1\linewidth]{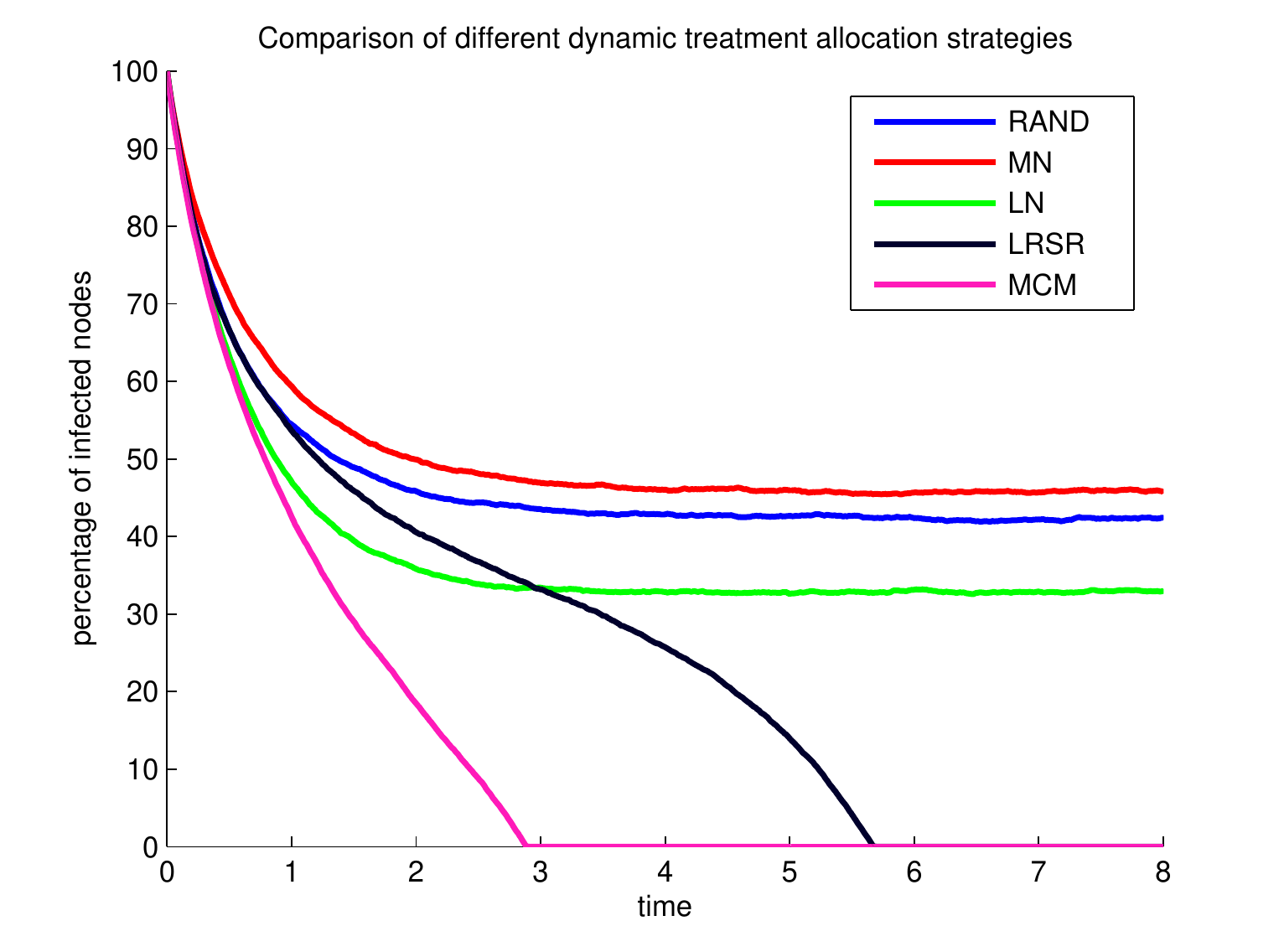}\vspace{-2mm}
		\caption{high resource efficiency: $e \op{=} 200$}
	\label{fig:twitter1}
	\end{subfigure}%
	\hspace{2em}
	\begin{subfigure}[b]{0.45\textwidth} 
		\centering
		\includegraphics[viewport=12 0 405 307, clip=true, width=1\linewidth]{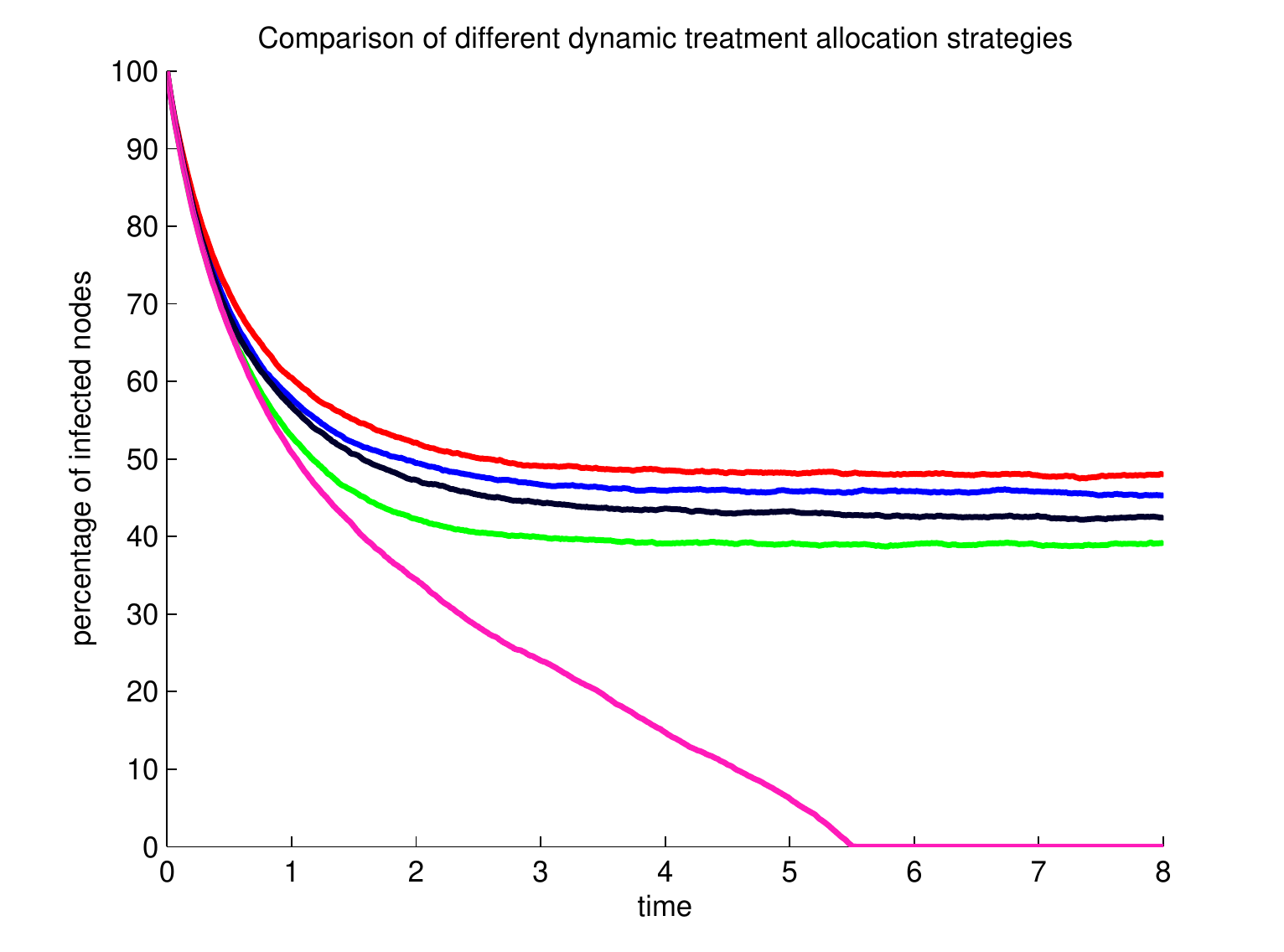}\vspace{-2mm}
	\caption{low resource efficiency: $e \op{=} 120$}
	\label{fig:twitter2}
	\end{subfigure}%
	\\
	\vspace{1em}
	\begin{subfigure}[b]{0.335\textwidth} 
		\centering
		\includegraphics[viewport=3 0 405 307, clip=true, width=1\linewidth]{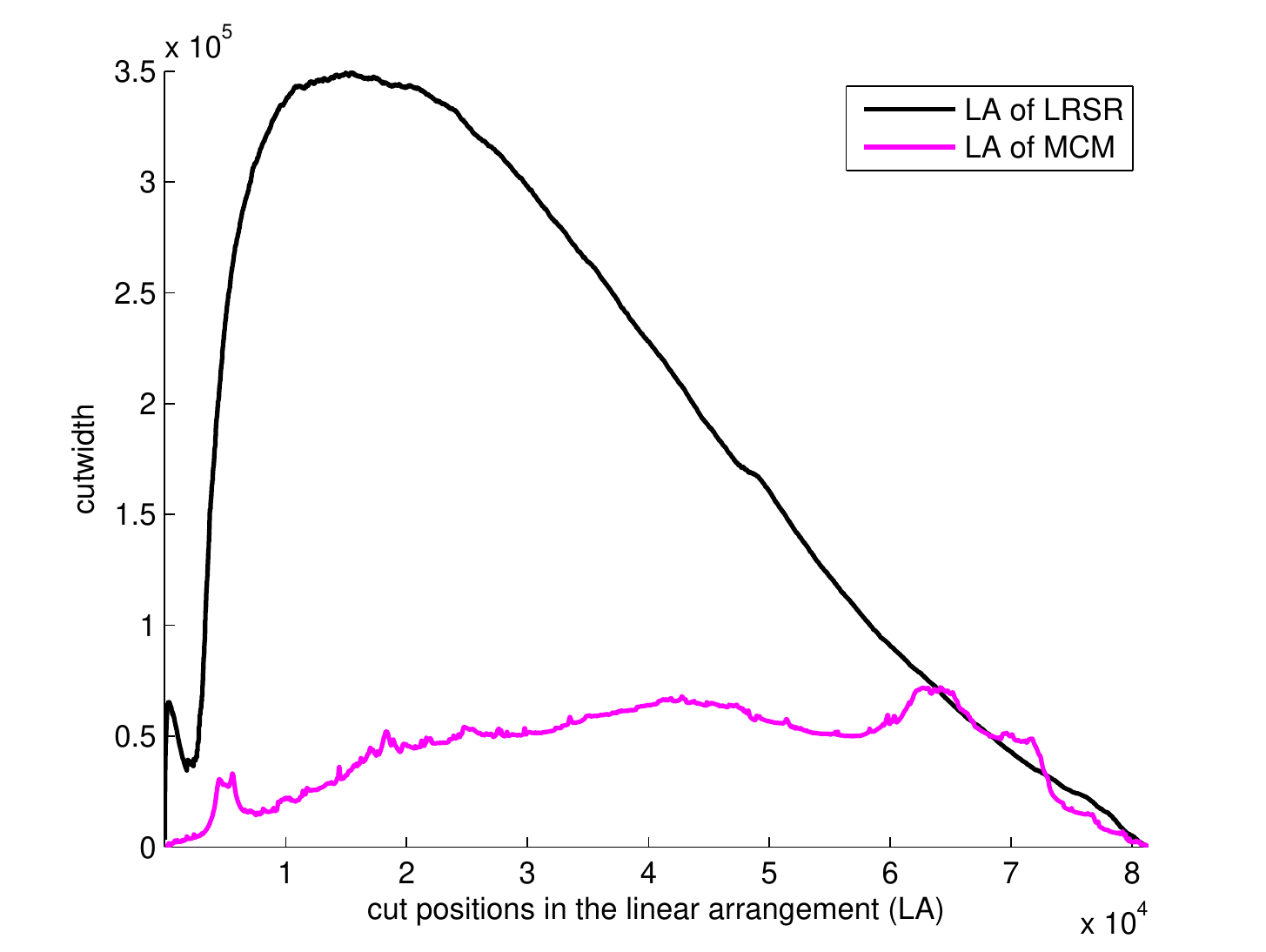}\vspace{-1mm}
		\caption{cutwidth values}
		\label{fig:cutwidths}
	\end{subfigure}
	\hspace{-0.2em}
	\begin{subfigure}[b]{0.31\textwidth} 
		\centering
		\includegraphics[viewport=119 273 470 570, clip=true, width=1\linewidth]{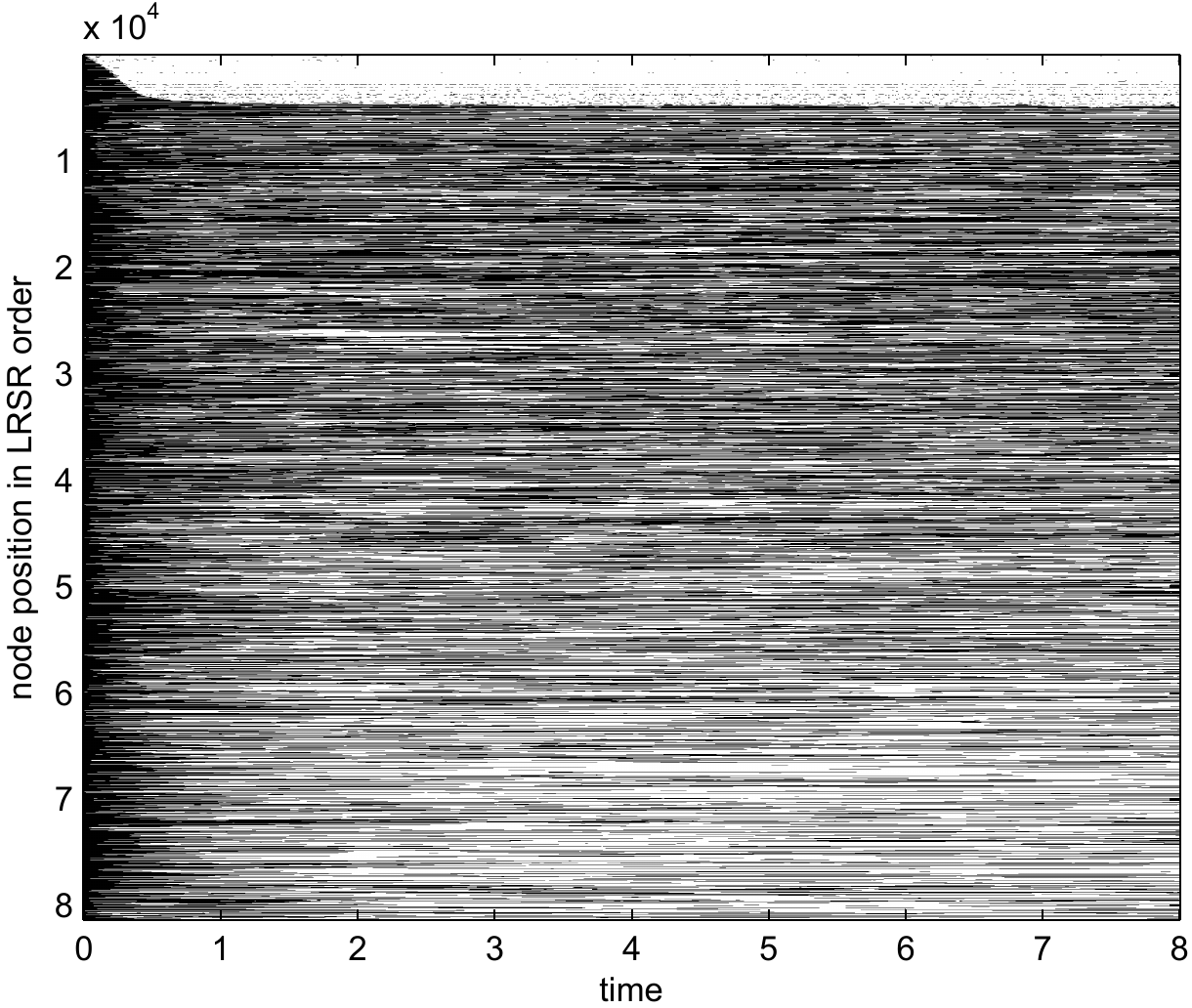}\vspace{-1mm}
		\caption{network state {--} LRSR}
		\label{fig:LRSRorder}
	\end{subfigure}
	\hspace{-0.2em}
	\begin{subfigure}[b]{0.31\textwidth} 
		\centering
		\begin{tikzpicture}
    \node[anchor=south west,inner sep=0] (image) at (0,0) {\includegraphics[viewport=119 273 470 570, clip=true, width=1\linewidth]{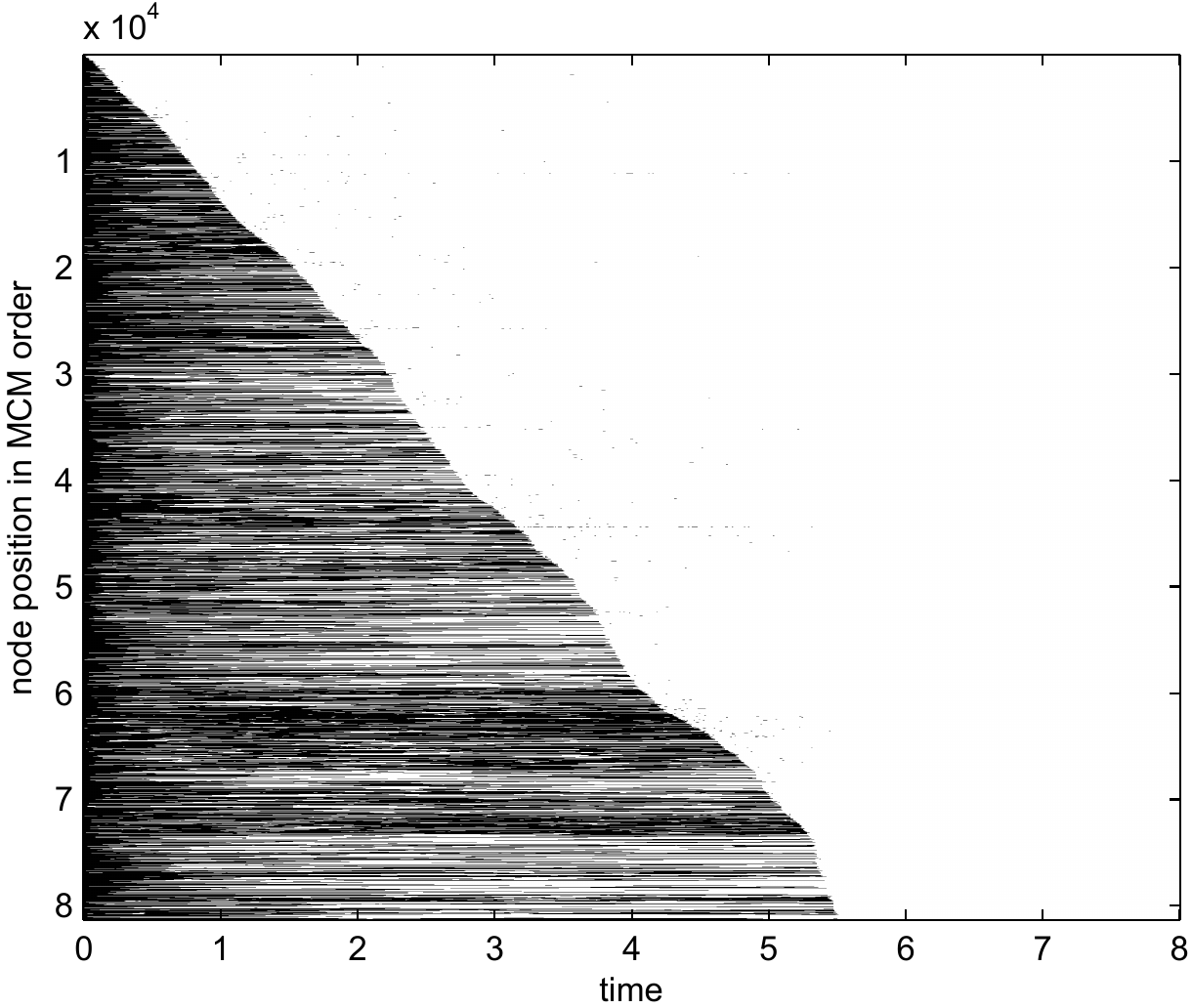}%
		\llap{\raisebox{5.5em}{\frame{\includegraphics[viewport=233 410 300 458, clip=true, width=0.5\linewidth]{orderMCM.pdf}}}}};
    \begin{scope}[x={(image.south east)},y={(image.north west)}]
        \draw[magenta,thick] (0.32,0.45) rectangle (0.519,0.605);
		\end{scope}
\end{tikzpicture}\vspace{-1mm}
		\caption{network state {--} MCM}
		\label{fig:MCMorder}
	\end{subfigure}
	\hspace{1em}
	\beforecaptvskip
	\beforecaptvskip
	\beforecaptvskip
	\caption{ 
	(\subref{fig:twitter1}-\subref{fig:twitter2})~Simulation of an undesired diffusion process in the TwitterNet subset of $81,306$ nodes \cite{conf/nips/McAuleyL12}. $\delta \op{=} 1$, $r \op{=} 0.1$ and $b_{tot} \op{=} 100$. MCM outperforms other heuristics. (\subref{fig:cutwidths})~Cutwidths for LRSR and MCM. (\subref{fig:LRSRorder}-\subref{fig:MCMorder})~Visualization of the diffusion in (\subref{fig:twitter2}) at the node level (contagious nodes in black). Nodes are ordered according to LRSR and MCM linear arrangements, respectively. A closer look of (\subref{fig:MCMorder}) is also provided in an inset figure.}
	\label{fig:expresults}
\end{figure}
For evaluating our strategy, the \emph{maximum cutwidth minimization} (MCM), we considered a realistic scenario where an SIS epidemic is being spread in a subset of the Twitter social network, denoted as TwitterNet 
\footnote{Available on: \texttt{http://snap.stanford.edu/data/index.html}.}. TwitterNet consists of $1,000$ ego-networks extracted from the social network \cite{conf/nips/McAuleyL12}. In order to perform experiments in the setting considered in this article, we symmetrized this network and obtained an undirected network of $81,306$ nodes and $1,342,303$ edges. The resulting network has a single connected component and contains a rich community structure.

The diffusion control literature is mainly focused on two types of strategies: static vaccination strategies and dynamic strategies based on a uniform mixing hypothesis. We compare MCM to the first type of strategies by considering state-of-the-art vaccination algorithms and use their output as a priority-order of the nodes, and to the second by considering random allocation of resources. As discussed in \Sec{sec:model}, the problem of finding the optimal budget through time is complementary to our approach. From this perspective, one of the primary question we address here is whether targeting specific nodes in the network can lead to substantial gains in the efficiency of the method, compared to the random distribution of the resources to the contagious nodes. \Fig{fig:twitter1} and \ref{fig:twitter2} compare the MCM strategy against the following four baseline approaches for creating a priority-order:
\begin{itemize}
\item \emph{Random order} (RAND): the order is a random permutation of the nodes of the network. 
\item \emph{Most neighbors} (MN): this strategy gives priority to high degree nodes. Intuitively, it begins by removing the contagion from the core of the network, and gradually reaches its periphery.
\item \emph{Least neighbors} (LN): this strategy gives priority to low degree nodes. Conversely to MN, LN begins with the periphery of the network and converges to its central part.
\item \emph{Least reduction in spectral radius} (LRSR): this strategy is based on the vaccination literature and more specifically the work in \cite{tong2012gelling}. This strategy gives priority to nodes whose removal will lead to the maximum decrease of the spectral radius of the adjacency matrix of the resulting network.
\end{itemize}

Concerning the applicability of \Theorem{th:boundExtinctionTime}, the maximum node degree in the network is $d_{max} \op{=} 3383$, which leads to a small $\epsilon$ value of $0.6$ even for the MCM strategy (and $\epsilon \op{=} 0.1$ for LRSR). We can thus be relatively confident in approximating the epidemic threshold under a given priority planning by $\frac{r\MaxCut}{b_{tot}}$. \Fig{fig:cutwidths} shows the cuts in the two best priority-orders, LRSR and MCM, while \Tab{tab:twittermaxcuts} summarizes the $\MaxCut$ values of the priority-orders produced by the different compared strategies. Note that the $\MaxCut$ of the priority-order produced by MCM is \emph{5 times smaller} than that of LRSR which is the best among its competitors. In essence, this implies that MCM would manage to suppress the diffusion process with resource efficiencies \emph{5 times smaller} than LRSR.

The results in \Fig{fig:twitter1} and \ref{fig:twitter2} illustrate that the proposed MCM strategy is more efficient than its competitors in removing the contagion from the network. In \Fig{fig:LRSRorder} and \ref{fig:MCMorder}, the evolution of the diffusion process is represented as follows: each line of the figure contains the state of one node of the network throughout the whole process (black for contagious and white for healthy). Nodes are sorted in y-axis according to the considered priority-order. We can observe that the cutwidth acts as a barrier for the LRSR: the large cutwidth values at the beginning of the LRSR order ($\frac{r\mathcal{C}_{max}}{b_{tot}} \op{\approx} e \op{=} 120$ for the $5000$-th node of the linear arrangement) prevents the strategy from consistently removing the contagion from more than the first $5000$ nodes of the plan. Note that healthy nodes also appear beyond the front, however this is not due to the control actions but rather due to self-recovery. 

On the contrary, MCM gradually reduces the contagion and the advancement of the front is clearly visible. These results agree with our previous analysis, and show that: i)~the uniform mixing hypothesis leads to a massive drop in efficiency, since MCM substantially outperforms the random allocation strategy, ii)~while efficient in the static vaccination problem, centrality-based priority-orders are suboptimal for the DRA problem, iii)~a good criterion for assessing the quality of a priority-order is actually in terms of its $\MaxCut$ value.
%
%
%
%=================================================================
\section{Conclusion}\label{sec:conclusion}
%=================================================================
%
In this paper, we presented a novel type of dynamic strategies for allocating resources in a network, called \emph{priority planning}, that aims to suppress an undesired contagion. We reduced the planning problem to that of linear arrangement of the nodes, and, based on theoretical analysis on the quality of any priority-order, i)~we demonstrated the key role of the \emph{maximum cutwidth} for assessing if a strategy would be eventually successful in removing the contagion, and ii)~we derived a strategy, called \emph{maximum cutwidth minimization} (MCM), that distributes resources to nodes in a priority-order with minimum maximum cutwidth. Our experimental results verified that, for a wide range of network types, the maximum cutwidth is indeed a good approximation of the epidemic threshold under a given strategy, and that the MCM strategy outperforms other competing strategies in a real-world social network.

\section*{Acknowledgments}
This research is part of the SODATECH project funded by the French Government within the program of ``\emph{Investments for the Future\,--\,Big Data}''.

\newpage
\bibliographystyle{ieeetr}
\bibliography{MCM2014}

\newpage
%
%
%==========================================
\section*{APPENDIX}
\subsection*{Mathematical arguments}
%==========================================
In order to simplify the demonstrations, and without loss of generality, we will reorder the nodes of $\mathcal{G}$ according to the priority-order $\lorder$, and thus consider that $\lorder(v_i) = i$.

%--------------------------------------------
\paragraph{Notations.}
%--------------------------------------------
%
First, let $X \in \{0,1\}^N$ be a state vector of size $N$(\ie describing the state of the network during the diffusion process), $\zero$ and $\one$ vectors of size $N$ that are all-zeros and all-ones, respectively, and $\bar{X} \op{=} \one \op{-} X$
. For $s \op{\subset} \{1,..., N\}$ a subset of nodes, let also $\one_s \op{=} (\one_{\{i \in s\}})_i$ be the indicator vector with ones for nodes in the set $s$. Then, we define $\tau_X$ as the extinction time of the contagion starting from the state $X$, i.e. the time needed for the Markov process to reach its absorbent state $X(t{=}\tau_X) \op{=} \zero$ when $X(t{=}0) \op{=} X$. We also denote the number of contagious nodes in network state $X$ as $N_I(X) \op{=} \one^\top X$, while $E_{I-S}(X) \op{=} X^\top A \bar{X}$ as the number of edges from a contagious to a healthy node, which edges we also refer to as \emph{contagious edges}.

The proof of \Theorem{th:boundExtinctionTime} 
relies on the following lemmas.

\begin{lemma}\label{th:increasingTau}
Under a priority plan, $X \mapsto \Exp{\tau_X}$ is monotonically increasing with respect to the natural partial order on $\{0, 1\}^N$ (\ie $X \leq Y$ if $X_i \leq Y_i$ for all $i$).
\end{lemma}
\vspace{-1em}
\begin{proof}
Let $X, Y \in \{0, 1\}^N$ be two initial states of the network such that $\forall i \in \{1,..., N\}, X_i \leq Y_i$. Let also $X(t)$ be the state vector of a contagion initially in state $X$, \ie $X(0) = X$, and $Y(t)$ be the state vector of a contagion initially in state $Y$, \ie $Y(0) = Y$. This lemma relies on the stochastic domination of $X(t)$ by $Y(t)$. This domination is due to the fact that, under the Markov Process defined by \Eq{eq:Xt}
, and when the control strategy is a priority planning, infection rates are increasing according to the natural partial order on $\{0, 1\}^N$, while recovery rates are decreasing. Thus, the initial inequality $X \leq Y$ will, in probability, grow during the contagion process.

The correct proof of this intuition relies on the strong monotonicity of the Markov Process, which we now prove. For all $X \leq Y$, the infection rate of each healthy node at state $X$ is lower than its infection rate at state $Y$, since $\sum_i A_{ji} X_i \leq \sum_i A_{ji} Y_i$. Also, the recovery rate of each infected node at state $X$ is higher than its recovery rate at state $Y$: if node $i$ is both infected in state $X$ and $Y$, and is treated in state $Y$, then $i = \min\{j\in\{1,...,N\}\!: Y_j = 1\}$ is the first infected node, and since the set of infected nodes of state $X$ is included in the set of infected nodes at state $Y$, we also have $i = \min\{j\in\{1,...,N\}\!: X_j = 1\}$ and $i$ is also treated at state $X$. We can thus apply Theorem 5.4 of \cite{masseyStoOrdering}, and $X(t)$ is a strongly monotone Markov process.

Let $X \op{\leq} Y$. If $X(t)$, $Y(t)$ are epidemic processes such that $X(0) \op{=} X$ and $Y(0) \op{=} Y$, then the strong monotonicity of the Markov process implies that $\forall t \op{\geq} 0, \mathbb{P}(\sum_i X_i(t) \op{=} 0) \op{\geq} \mathbb{P}(\sum_i Y_i(t) \op{=} 0)$, which may be rewritten as follows: $\mathbb{P}(\tau_X \op{\leq} t) \op{\geq} \mathbb{P}(\tau_Y \op{\leq} t)$. This means that $\tau_Y$ dominates $\tau_X$ and thus $\Exp{\tau_X} \op{\leq} \Exp{\tau_Y}$.
\end{proof}

\begin{lemma}\label{th:generalBound}
Assume $b_{tot} = 1$ and let $X_n^j$ the worst state vector after $j$ additional infections from the $n^{th}$ state of the planned strategy $X_n$:
\begin{equation}
X_n^j = \argmax_{
\begin{array}{l}
\one_s ~\st\\
\{n,..., N\} \subset s,\\
|s \cup \{1,..., n-1\}| = j
\end{array}}
{\Exp{\tau_{\one_s}}}.
\end{equation}
Then the following bound for the expected extinction time under the priority planning and starting from a total infection holds:
\\
$\forall K \geq 1$ and $\rho > \beta \left[\sum_{n=1}^N \prod_{j=0}^K E_{I-S}(X_n^j)\right]^{\frac{1}{K+1}} - \delta$,
\begin{equation}
\Exp{\tau_\one} \ \leq \ \frac{\sum_{k=0}^K f(k)}{(\rho + \delta)(1 - f(K+1))},
\end{equation}
where
\begin{equation}
f(k) = \sum_{n=1}^N \left(\frac{\beta}{\rho + \delta}\right)^k \prod_{j=0}^{k-1} E_{I-S}(X_n^j).
\end{equation}
\end{lemma}
\vspace{-1em}
\begin{proof}
For every state vector $X$, we have:
\begin{equation}
\Exp{\tau_X} = \Exp{t_1 + \tau_{X(t_1)}},
\end{equation}
where $t_1$ is the time of the first change in the state vector. Three types of events can happen: either a node recovers by itself (at a rate $\delta$), or a node is healed by a resource (at a rate $\rho$ \op{+} $\delta$), or a node is infected (at a rate $\beta$). The number of infected nodes is $N_I(X)$. The first contagious node, denoted as $i_X \op{=} \min\{j\in\{1,...,N\}\!: X_j \op{=} 1\}$, receives a resource, and the number of nodes that can be infected is $E_{I-S}(X)$. Thus,
\begin{equation}
\begin{array}{ll}
\Exp{\tau_X} = &\frac{1}{\delta N_I(X) + \rho + \beta E_{I-S}(X)}\\\\
               &+ \frac{\delta (N_I(X) - 1)}{\delta N_I(X) + \rho + \beta E_{I-S}(X)}\Exp{\tau_{X(t_1)} | \mbox{self-recovery of a node at } t_1}\\\\
               &+ \frac{\rho + \delta}{\delta N_I(X) + \rho + \beta E_{I-S}(X)}\Exp{\tau_{X - \one_{i_X}}}\\\\
               &+ \frac{\beta E_{I-S}(X)}{\delta N_I(X) + \rho + \beta E_{I-S}(X)}\Exp{\tau_{X(t_1)} | \mbox{infection at } t_1}.\\
\end{array}
\end{equation}

Using \Lemma{th:increasingTau}, we get that $\Exp{\tau_{X(t_1)} | \mbox{self-recovery of a node at } t_1} \leq \Exp{\tau_X}$, which leads to:
\begin{equation}
(\delta + \rho + \beta E_{I-S}(X))\Exp{\tau_X} \leq 1 + (\rho + \delta) \Exp{\tau_{X - \one_{i_X}}} + \beta E_{I-S}(X) \Exp{\tau_{X(t_1)} | \mbox{infection at } t_1}.
\end{equation}

Let $u_n^j \op{=} \Exp{\tau_{X_n^j}}$. For all $j \op{\geq} 1$, by definition of $X_n^{j-1}$ and $X_n^{j+1}$, we have $\Exp{\tau_{X_n^j \op{-}\one_{i_{X_n^j}}}} \op{\leq} u_n^{j-1}$ and $\Exp{\tau_{X_n^j(t_1)} | \mbox{infection at } t_1} \op{\leq} u_n^{j+1}$. This comes from the fact that, as the order is static, the $j$ infected nodes that are among $\{1,..., n\op{-}1\}$ will receive a resource first. We thus get the following recurrence equation for the $u_n^j$:
\begin{equation}
(\delta + \rho)(u_n^j - u_n^{j-1}) \leq 1 + \beta E_{I-S}(X_n^j)(u_n^{j+1} - u_n^j),
\end{equation}
which can be iterated:
\begin{equation}
\begin{array}{lll}
(\rho + \delta) (u_n^0 - u_{n+1}^0) &\leq &\sum_{k=0}^{K-1} (\frac{\beta}{\rho + \delta})^k \prod_{j=0}^{k-1} E_{I-S}(X_n^j)\\\\
                                         &&+ (\rho + \delta) (u_n^{K+1} - u_n^K) (\frac{\beta}{\rho + \delta})^K \prod_{j=0}^K E_{I-S}(X_n^j)\\\\
                                    &\leq &\sum_{k=0}^{K-1} (\frac{\beta}{\rho + \delta})^k \prod_{j=0}^{k-1} E_{I-S}(X_n^j)\\\\
																		     &&+ (\rho + \delta) u_1^0 (\frac{\beta}{\rho + \delta})^K \prod_{j=0}^K E_{I-S}(X_n^j),\\
\end{array}
\end{equation}
since $u_n^{K+1} \op{\leq} u_1^0$ using \Lemma{th:increasingTau} and $u_1^0 \op{=} \Exp{\tau_\one}$.

We can now derive the final formula by summing over $n$ and using the definition $f(k) = \sum_{n=1}^N (\frac{\beta}{\rho + \delta})^k \prod_{j=0}^{k-1} E_{I-S}(X_n^j)$:
\begin{equation}
(\rho + \delta) (1 - f(K+1)) \Exp{\tau_\one} \leq \sum_{k=0}^{K-1} f(k).
\end{equation}
\end{proof}
\begin{lemma}\label{th:EisBound}
For $n \op{\in} \{1,..., N\}$ and $j \op{\in} \{0,..., n\op{-}1\}$,
\begin{equation}
E_{I-S}(X_n^j) \ \leq \ E_{I-S}(X_n) + j d_{max},
\end{equation}
where $d_{max} \op{=} \max_i \sum_j A_{ij}$ is the highest degree of the network.
\end{lemma}
\vspace{-1em}
\begin{proof}
The contagious nodes of $X_n^j$ are the contagious nodes of $X_n$ and exactly $j$ additional nodes. Since a node can have at most $d_{max}$ neighbors, then each of the $j$ additional nodes can add at most $d_{max}$ edges to the set of contagious edges of the network.
\end{proof}

\begin{lemma}\label{th:xiBound}
Let $a \op{\geq} 0$ and $\xi$ be the (unique) positive solution to $\xi \op{-} \ln(1+\xi) \op{=} a$. The following inequality holds:
\begin{equation}
\xi \leq a + 2 \sqrt{a}.
\end{equation}
\end{lemma}
\vspace{-1em}
\begin{proof}
$x \op{-} \ln(1\op{+}x)$ is convex, thus always above its tangent line:
\\
$\forall x_0 > 0,$ 
\begin{equation}
a = \xi - \ln(1+\xi) \geq (x_0 - \ln(1+x_0)) + \frac{x_0}{1+x_0}(\xi - x_0),
\end{equation}
and thus,
\begin{equation}
\begin{array}{ll}
\xi &\leq \frac{1+x_0}{x_0}(a + \ln(1+x_0)) - 1\\\\
    &\leq \frac{1+x_0}{x_0} a + x_0.\\
\end{array}
\end{equation}
The final result is obtained by setting $x_0 \op{=} \sqrt{a}$.
\end{proof}

We can now prove the \Theorem{th:boundExtinctionTime} 
using the above lemmas.
\begin{proof}[\textbf{Proof of \Theorem{th:boundExtinctionTime}}]
Using \Lemma{th:generalBound} and \Lemma{th:EisBound}, we obtain a bound on the extinction time depending on $\MaxCut \op{=} \max_n E_{I-S}(X_n)$,
\\
$\forall K \geq 1$ and $\rho > \beta \left[N \prod_{j=0}^K (\MaxCut + j d_{max})\right]^{\frac{1}{K+1}} - \delta$,
\begin{equation}
\Exp{\tau_\one} \ \leq \ \frac{\sum_{k=0}^K f(k)}{(\rho + \delta)(1 - f(K+1))}
\end{equation}
where
\begin{equation}
f(k) = N (\frac{\beta}{\rho + \delta})^k \prod_{j=0}^{k-1} (\MaxCut + j d_{max}),
\end{equation}
using the approximation $\sum_{n=1}^N \prod_{j=0}^{k-1} (E_{I-S}(X_n) + j d_{max}) \leq N \prod_{j=0}^{k-1} (\MaxCut + j d_{max})$.

Finally, we need to select a proper value for $K$ and derive the final result. Let $\xi$ be the unique solution of $\xi \op{-} \ln(1\op{+}\xi) \op{=} \frac{d_{max} \ln{N}}{\MaxCut}$ and $K^* \op{=} \lfloor \frac{\MaxCut}{d_{max}} \xi\rfloor$. Using the particular value of $K^*$,
\begin{equation} \label{eq:justaninequality}
\begin{array}{ll}
\sum_{j=0}^{K^*} \ln(1 + j \frac{d_{max}}{\MaxCut})
    &\leq \int_{0}^{K^* + 1} \ln(1 + x \frac{d_{max}}{\MaxCut})dx\\
    &= (K^* + 1 + \frac{\MaxCut}{d_{max}})\ln(1 + (K^* + 1) \frac{d_{max}}{\MaxCut}) - (K^* + 1)\\
		&\leq (K^* + 1) \ln(1 + (K^* + 1) \frac{d_{max}}{\MaxCut}) + \frac{\MaxCut}{d_{max}}(\ln(1 + \xi) - \xi)\\
    &= (K^* + 1) \ln(1 + (K^* + 1) \frac{d_{max}}{\MaxCut}) - \ln(N),\\
\end{array}
\end{equation}
where the second inequality is due to $\frac{\MaxCut}{d_{max}}(K^* \op{+} 1) \op{\geq} \xi$ and the monotonic decrease of $x \op{\mapsto} \ln(1\op{+}x) \op{-} x$ for $x \op{\geq} 0$.
\\
From \Eq{eq:justaninequality}, we derive that $f(K^*\op{+}1) \op{\leq} \left[\frac{\beta}{\rho + \delta}(\MaxCut \op{+} (K^*\op{+}1) d_{max})\right]^{K^*+1}$. We thus have:
\\
For $\rho \op{>} \beta (\MaxCut \op{+} (K^*+1) d_{max}) \op{-} \delta$,
\begin{equation}
\begin{array}{ll}
\Exp{\tau_\one} &\leq \frac{\sum_{k=0}^{K^*} f(k)}{(\rho + \delta)(1 - f(K^*+1))}\\\\
                &\leq \frac{N\sum_{k=0}^{K^*} \left[\frac{\beta}{\rho + \delta} (\MaxCut + (K^*+1) d_{max})\right]^k}{(\rho + \delta)(1 - f(K^*+1))}\\\\
                &\leq \frac{N}{(\rho + \delta)(1 - \frac{\beta}{\rho + \delta} (\MaxCut + (K^*+1) d_{max}))}\cdot\frac{1 - \left[\frac{\beta}{\rho + \delta} (\MaxCut + (K^*+1) d_{max})\right]^{K^*+1}}{1 - f(K^*+1)}\\\\
                &\leq \frac{N}{\rho + \delta - \beta (\MaxCut + (K^*+1) d_{max})}.\\
\end{array}
\end{equation}
\\
Finally, using \Lemma{th:xiBound}, $d_{max} K^* \op{\leq} \MaxCut \xi \op{\leq} \MaxCut \left[\frac{d_{max} \ln{N}}{\MaxCut} \op{+} 2\sqrt{\frac{d_{max} \ln{N}}{\MaxCut}}\right]$,\\ which proves the desired bound.
\end{proof}

\end{document}

%% file: definitions_arxiv.tex
\usepackage{graphics}

\usepackage{algorithm}
\usepackage{algorithmic}
\usepackage{enumitem}

\usepackage{amsthm}
\newtheorem{theorem}{Theorem}
\newtheorem{lemma}{Lemma}
\newtheorem{corollary}{Corollary}

\usepackage{hyperref}
\usepackage{url}
\usepackage{amsmath}

\usepackage{amsfonts}
\usepackage{dsfont}
\usepackage{upgreek}
\usepackage{bbm}			% mathbb vs mathds
\usepackage{mathbbol}

\usepackage{graphicx}
\usepackage{epstopdf}
\usepackage{tikz}
\usepackage{float}

\usepackage{caption}
\usepackage{subcaption}

\usepackage[square, sort, numbers]{natbib}
	
\usepackage{xspace}

\usepackage[scaled=.7]{beramono}

\newcommand{\Sec}[1]		{Sec.\,\ref{#1}}
\newcommand{\Fig}[1]		{Fig.\,\ref{#1}}

\newcommand{\Eq}[1]			{Eq.\,\ref{#1}}
\newcommand{\Tab}[1]		{Tab.\,\ref{#1}}
\newcommand{\Alg}[1]		{Alg.\,\ref{#1}}
\newcommand{\Theorem}[1]{Theorem~\ref{#1}}

\newcommand{\Corollary}[1]{Corollary~\ref{#1}}
\newcommand{\Lemma}[1]{Lemma~\ref{#1}}
\newcommand{\ie}   			{i.e.\ }
\newcommand{\eg}   			{e.g.\ }
\newcommand{\wrt}   		{w.r.t.\ }
\newcommand{\st}   			{\mbox{s.t.\ }}
\newcommand{\Erdos}   	{Erd\"os-R\'enyi }
\newcommand{\one}       {\mathbb{1}}
\newcommand{\zero}      {\mathbb{0}}

\newcommand{\Exp}[1]    {\mathbbm{E}[#1]}
\newcommand{\real}      {\mathbbm{R}}
\newcommand{\op}[1]     {\,{#1}\,}

\newcommand{\spinlinetitle}[2]  {\vspace{0.5em}\noindent\textbf{\emph{#1}{#2}}}

\newcommand{\beforecaptvskip} {\vskip -0.07in}

\newcommand{\formulastyle}   	{\textstyle} 

\DeclareMathOperator*{\argmax}{\mbox{argmax}}

\newcommand{\lorder}	{\ell}
\newcommand{\objf}  	{\phi}

\newcommand{\MaxCut}  {\mathcal{C}_{max}}